\documentclass[12pt]{article}
\usepackage{amsmath,amsfonts,amsthm,amssymb}
\usepackage[hscale=0.7, vscale=0.76]{geometry}
\usepackage{graphicx,color}
\newtheorem{thm}{Theorem}
\newtheorem{cor}{Corollary}
\theoremstyle{remark}
\newtheorem{rem}{Remark}
\let\bx\square
\def\R{\mathbb R}
\def\D{\mathrm D}
\def\d{\mathrm d}
\begin{document}
\title{On the existence of the M{\o}ller wave operator for wave equations with small dissipative terms.}
\author{Jens Wirth}
\maketitle
\begin{abstract}
The aim of this short note is to reconsider and to extend a former result of K.~Mochizuki \cite{Moc76}, \cite{MN96}
on the existence of the scattering operator for wave equations with small dissipative terms.

Contrary to the approach used by Mochizuki we construct the wave operator explicitly in terms of the parametrix construction
obtained by a (simplified) diagonalization procedure, cf. \cite{Yag97}. The method is based on ODE techniques.

These considerations are part of a larger project and the idea is taken from \cite{Wir02B} and 
generalized to $x$-dependent coefficients.
\vspace{1cm}

{\sl \noindent
keywords: Cauchy problem, wave equation, scattering operator, dissipative term}
\end{abstract}

We consider the Cauchy problem 
\begin{equation}\label{eq:CP}
 \bx u+b(t,x)u_t=0,\quad u(0,\cdot)=u_1,\quad \D_t u(0,\cdot)=u_2 
\end{equation}
with $b\in L_1(\R,L_\infty(\R^n))\cap L_\infty(\R^{1+n})$. We restrict our 
calculations to spaces of dimension $n\ge2$. The modification to obtain results also for $n=1$ are 
obvious.

We denote by $E=\dot H^1\times L_2$ the energy space. We prove that in the 
energy space $(u,\D_t u)$ converges to the local energy of a 
solution $\tilde u$ of the free wave equation $\bx \tilde u=0$. As usual we 
use the representation $\dot H^1=|\D|^{-1}L_2$, $|\D|^{-1}$ beeing the Riesz potential
operator of order $1$, inverse to the operator $|\D|=\sqrt{-\Delta}$. In our calculations we use this isomorphism $\dot H^1\simeq L_2$ to restrict ourselves to calculations in $L_2$-space.

\begin{thm}\label{thm1}
  Assume $b\in L_1(\R,L_\infty(\R^n))\cap L_\infty(\R^{1+n})$.

  There exist isomorphisms $W_\pm:E\to E$ of the energy space such that for
  $u=u(t,x)$ the solution of \eqref{eq:CP} to data $(u_1,u_2)\in E$ and
  for $(\tilde u_1,\tilde u_2)=W_\pm(u_1,u_2)$ and $\tilde u$ the solution of
  the free wave equation $\bx\tilde u=0$ to data $\tilde u(0,\cdot)=\tilde u_1$,   
  $\D_t\tilde u(0,\cdot)=\tilde u_2$ the asymptotic relation
  $$ ||(u,\D_tu)-(\tilde u,\D_t\tilde u)||_E\to 0 
  \qquad\text{ as }\quad t\to\pm\infty $$
  holds.  
\end{thm}

We subdivide the proof in some steps and construct these wave operators explicitly in terms of the solution representation.

Let $U=(|\D|\hat u,\D_t\hat u)^T$. Then $U$ satisfies the equation
$$ \D_tU=\begin{pmatrix}&|\D|\\|\D|&\end{pmatrix}U
   +\begin{pmatrix}\quad &\\&ib(t,x)\end{pmatrix}U. $$

The first matrix operator maps $\dot H^{1}\to L_2$ while due to our assumptions the second one maps
$L_2\to L_2$. We will understand the first operator as closed unbounded operator on $L_2$ with domain 
$\dot H^1$ and solve in a first step the corresponding evolution equation. 
In a second step we understand the second operator as small perturbation of the first one.

We diagonalize the first matrix operator. Therefore we use
$$ M=\begin{pmatrix} 1&-1\\1&1\end{pmatrix},
  \qquad M^{-1}=\frac12\begin{pmatrix}1&1\\-1&1\end{pmatrix} $$
and consider $U^{(0)}=M^{-1}U$. We get
$$ \D_tU^{(0)}=M^{-1}\begin{pmatrix}&|\D|\\|\D|&\end{pmatrix}MU^{(0)}
   +M^{-1}\begin{pmatrix}\quad &\\&ib(t,x)\end{pmatrix}M U^{(0)}=\mathcal D U^{(0)}+B(t,x)U^{(0)}, $$
where $\cal D$ is the diagonal operator
$$ \mathcal D =\begin{pmatrix}|\D|&\\&-|\D|\end{pmatrix} $$
and $B(t,x)\in L_1(\R,L_\infty(\R^n))$ is a matrix. Multiplication by this matrix defines a bounded operator
on $L_2$ with $||B(t,x)\cdot||_{2\to2}=||B(t,\cdot)||_\infty\in L_1(\R)$.

We start by solving the operator valued evolution equation
$$\D_t\mathcal E_0(t,s)=\mathcal D\mathcal E_0(t,s),\qquad \mathcal E_0(s,s)=I:L_2\to L_2.$$
We look for a solution from the evolution space
$$ \mathcal E_0(\cdot,s)\in C(\R,L_2\to L_2)\cap C^1(\R,L_2\to\dot H^{-1}) .$$

This solution is given by the fundamental solution corresponding to the wave equation
$$ \mathcal E_0(t,s)=e^{i(t-s)\mathcal D}
  =\begin{pmatrix}e^{i(t-s)|D|}&\\&e^{-i(t-s)|D|}\end{pmatrix}. $$
The operator $\mathcal E_0$ is unitary as well as $U_0(t,s)$, obtained
by a similarity transform of $M\mathcal E_0(t,s)M^{-1}$ with the (canonical) isomorphism $E\simeq L_2$.
The operator $U_0(t)$ is the unitary solution operator of the homogeneous wave equation in the 
energy space.

In a second step we construct the solution to $\D_t-\mathcal D-B(t,x)$. 
Therefore let
$$ \mathcal R(t,s)=\mathcal E_0(s,t)B(t,x)\mathcal E_0(t,s)$$
(as concatenation of (bounded) operators on $L_2$)
and 
$$
 \mathcal Q(t,s)=I+\sum_{k=1}^\infty i^k \int_s^t\mathcal R(t_1,s)
  \int_s^{t_1}\mathcal R(t_2,s)\dots\int_s^{t_{k-1}}\mathcal R(t_k,s)\d t_k\dots\d t_1 
$$
in the sense of Bochner integrals.
The matrix operator $\mathcal Q(t,s)$ solves the Cauchy problem
$$ \D_t \mathcal Q(t,s)-\mathcal R(t,s)\mathcal Q(t,s)=0,\qquad \mathcal Q(s,s)=I:L_2\to L_2. $$

Using $\mathcal Q(t,s)$ we can express the fundamental solution to the
diagonalized system. Let $\mathcal E(t,s)=\mathcal E_0(t,s)\mathcal Q(t,s)$. Then we obtain
\begin{align*}
  \D_t (\mathcal E_0\mathcal Q)&
  =(\D_t\mathcal E_0)\mathcal Q+\mathcal E_0(\D_t\mathcal Q)
  =\mathcal D\mathcal E_0\mathcal Q+\mathcal E_0\mathcal R(t,s)\mathcal Q\\
  &=\mathcal D\mathcal E_0\mathcal Q+B(t,x)\mathcal E_0\mathcal Q
\end{align*}
and $\mathcal E_0(s,s)\mathcal Q(s,s)=I$. Thus $\mathcal E(t,s)$ is the desired 
fundamental solution.

Hence $M\mathcal E(t,s)M^{-1}$ is related to the operator
$$ U(t,s):E\ni(u(s),\D_t u(s))\mapsto(u(t),\D_t u(t))\in E$$
for solutions $u$ to $\bx u+b(t,x)u_t=0$.

We estimate the norm of this operator. We do this step by step. At first we have
$$ ||\mathcal E_0(t,s)||=1. $$
The next estimate is
$$ ||\mathcal R(t,s)||\le ||B(t,\cdot)||_\infty \in L_1(\R) $$
which will be used to estimate $\mathcal Q(t,s)$. We use the following statement
$$ 
  \left|\int_s^t r(t_1,s) \int_s^{t_1}r(t_s,s)\dots\int_s^{t_{k-1}}r(t_k,s)\d t_k\dots\d t_1 \right|\le \frac1{k!}\left(\int_s^t|r(\tau,s)|\d\tau\right)^k.
$$
This can be proved using induction over $k$. Combined with the series representation of $\mathcal Q$ we get
\begin{align*}
 ||\mathcal Q(t,s)-I||
 &\le\sum_{k=1}^\infty\frac1{k!}\left(\int_s^t||B(\tau,\cdot)||_\infty\d\tau\right)^k\\
 &=\exp\left\{\int_s^t||B(\tau,\cdot)||_\infty\d\tau\right\}-1\le C
\end{align*}
and therefore 
$$ ||\mathcal E(t,s)||\le C. $$

Integrability of $||B(t,\cdot)||_\infty$ implies further $\mathcal Q(t,s)\to I$ as $t,s\to\infty$.

\begin{rem}
  The constant in this estimate can be larger than $1$. This is due to the fact that we have not
  required the condition $b(t,x)\ge0$.
\end{rem}

\begin{rem}
We are interested in the M{\o}ller wave operator $W_+$. This operator can be understood as a limit in
the following sense. We consider data $(u_1,u_2)$ from the energy space and apply the solution
operator $U(t,0)$. Then we go back to the initial line using the solution operator of the homogeneous problem $U_0(-t)$. This gives data to the homogeneous wave equation which produce a 
solution coinciding with $u$ at the time level $t$. Now we let $t\to\infty$
$$W_+=\lim_{t\to\infty} U_0(-t)U(t,s).$$
If this limit exists in some sense, we have constructed the first M{\o}ller wave operator.
The second one, $W_-$ will be obtained by replacing $t$ by $-t$. 
In our case we will see that these limits exists as strong limits in the operator norm.
\end{rem}

From $\mathcal E_0(0,t)\mathcal E(t,0)=\mathcal Q(t,0)$ it seems natural to ask whether $\lim_{t\to\infty} \mathcal Q(t,0)$ exists in $L_2\to L_2$. Therefore we
consider the difference
\begin{align*}
  \mathcal Q(t,0)-\mathcal Q(s,0)
  &=\sum_{k=1}^\infty i^k\bigg[\int_0^t\mathcal R(t_1,0)
  \int_0^{t_1}\mathcal R(t_2,0)\dots\int_0^{t_{k-1}}\mathcal R(t_k,0)\d t_k\dots\d t_1 \\
   &\qquad\qquad -\int_0^s\mathcal R(t_1,0,\xi)
  \int_0^{t_1}\mathcal R(t_2,0)\dots\int_0^{t_{k-1}}\mathcal R(t_k,0)\d t_k\dots\d t_1
  \bigg]\\
  &=\sum_{k=1}^\infty i^k\int_s^t\mathcal R(t_1,0)
  \int_0^{t_1}\mathcal R(t_2,0)\dots\int_0^{t_{k-1}}\mathcal R(t_k,0)\d t_k\dots\d t_1.
\end{align*}
If we apply $||\cdot||_{2\to2}$ on both sides and use the same statement as above to estimate the integrals we get
\begin{align*}
  ||\mathcal Q(t,0,\cdot)-\mathcal Q(s,0,\cdot)||_{2\to2}
  &\le\sum_{k=1}^\infty \int_s^t ||B(t_1,\cdot)||_\infty
   \frac1{(k-1)!}\left(\int_0^{t_1}||B(\tau,\cdot)||_\infty\d\tau\right)^{k-1}\d t_1\\
  &\le\int_s^t||B(t_1,\cdot)||_\infty\sum_{k=0}^\infty \frac1{k!}\left(\int_0^{t_1}||B(\tau,\cdot)||_\infty\d\tau\right)^k\d t_1\\
  &=\int_s^t||B(t_1,\cdot)||_\infty\exp\left\{\int_0^{t_1}||B(\tau,\cdot)||_\infty\d\tau\right\}\d t_1\to0
\end{align*}
as $t,s\to\infty$ from the integrability of $||B(t,\cdot)||_\infty$.

Thus the limit exists in $L_2\to L_2$ and we can define
$$ \tilde W_+=\lim_{t\to\infty} M\mathcal Q(t,0)M^{-1}. $$

\begin{rem}
By the construction it follows that
$$ W_+=\lim_{t\to\infty} U_0(-t)U(t,0) $$
is related to $\tilde W_+$ via the isomorphism $E\simeq L_2$.
\end{rem}

The transpose\footnote{In the usual matrix sense.} of the inverse of $Q(t,s)$ satisfies the related equation
$$ \D_t-\mathcal Q^{-T}(t,s)+\mathcal R^T(t,s)\mathcal Q^{-T}(t,s)=0,\qquad \mathcal Q^{-T}(s,s)=I. $$
Thus we can estimate $\mathcal Q^{-T}$ in a similar style as $\mathcal Q$, especially we can prove that
$$ \lim_{t\to\infty} \mathcal Q^{-1}(t,s) $$
exists.

The same argumentation is vaild for $t\to-\infty$ and defines the second operator $W_-$. Especially the scattering operator $S=W_+W_-^{-1}$ exists.

\begin{cor}
  Under the assumptions of Theorem \ref{thm1} exists the scattering operator $S:E\to E$ and is invertible.
\end{cor}

\begin{figure}[h]
\begin{center}
\input{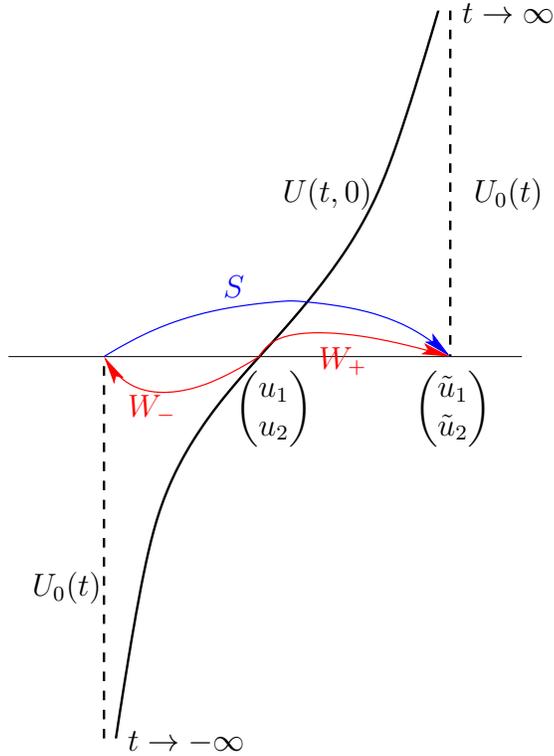}
\end{center}
\caption{Short overview on operators involved in this note.}
\end{figure}

We can also estimate the rate of convergence in Theorem \ref{thm1}.
\begin{cor}
  Under the assumptions of Theorem \ref{thm1} it holds
  \[
  ||(u,\D_tu)-(\tilde u,\D_t\tilde u)||_E\lesssim 
  ||(u_1,u_2)||_E\;\int_t^\infty ||b(\tau,\cdot)||_\infty\d\tau.
  \]
\end{cor}
\begin{proof}
  The statement follows directly from
  $$
    \mathcal Q(\infty,0)-\mathcal Q(t,0)=\sum_{k=1}^\infty i^k\int_t^\infty\mathcal R(t_1,0)
  \int_0^{t_1}\mathcal R(t_2,0)\dots\int_0^{t_{k-1}}\mathcal R(t_k,0)\d t_k\dots\d t_1.
  $$
  and
  \begin{align*}
   ||\mathcal Q(\infty,0)-\mathcal Q(t,0)||_{2\to2}
   &\le \int_t^\infty ||B(t_1,\cdot)||_\infty
       \exp\left\{\int_0^{t_1}||B(\tau,\cdot)||_\infty\d\tau\right\}\d t_1\\
   &\lesssim \int_t^\infty ||b(\tau,\cdot)||_\infty\d\tau,
   \end{align*}
  where $\mathcal Q(\infty,s)=\lim_{t\to\infty}\mathcal Q(t,s)$.
\end{proof}

For our considerations it was essential that we had \emph{not} to indroduce a subdivision of the phase space into zones. This enables
us to give a definition of the wave operator $W_+(\xi)$ globally in the phase variable $\xi$. 

For a more general treatment including a subdivision of the phase space we refer to \cite{Wir02B} and the modified scattering result
discussed there.

\end{document}